\documentclass[12pt, times]{article}
\usepackage{graphicx}
\usepackage{amssymb,amsmath}
\def\G{\mathcal{G}}

\def\R{\mathcal{R}}

\def\T{\mathcal{T}}
\def\R{\mathcal{R}}
\def\N{\mathcal{N}}
\def\M{\mathcal{M}}
\def\Aut{\mathrm{Aut}}
\newtheorem{defn}{Definition}
\newtheorem{thm}{Theorem}
\newtheorem{corol}{Corollary}
\newtheorem{lemma}{Lemma}

\newtheorem{propos}{Proposition}

\def\proof{\noindent{\it Proof}\/.\ }

\begin{document}
\title{\bf Characterization and Enumeration of Toroidal $K_{3,3}$-Subdivision-Free Graphs}
\bigbreak
\author{Andrei Gagarin, Gilbert Labelle and Pierre Leroux \\[0.1in]
\small \it Laboratoire de Combinatoire et d'Informatique Math\'ematique (LaCIM),\\ 
\small \it Universit\'e du Qu\'ebec \`a Montr\'eal, Montr\'eal, Qu\'ebec, CANADA, H3C 3P8\\[0.1in]
\small e-mail: {\it gagarin@math.uqam.ca, labelle.gilbert@uqam.ca} and \it leroux.pierre@uqam.ca }
\maketitle
\begin{abstract}
We describe the structure of $2$-connected non-planar toroidal graphs with no $K_{3,3}$-subdivisions, using an appropriate substitution of planar networks into the edges of certain graphs called toroidal cores. The structural result is based on a refinement of the algorithmic results for graphs containing a fixed $K_5$-subdivision in [A. Gagarin and W. Kocay, ``Embedding graphs containing $K_5$-sub\-divisions'', Ars Combin. {\bf 64} (2002), 33-49]. It allows to recognize these graphs in linear-time and makes possible to enumerate labelled $2$-connected toroidal graphs containing no $K_{3,3}$-sub\-divisions and having minimum vertex degree two or three by using an approach similar to [A. Gagarin, G. Labelle, and P. Leroux, "Counting labelled projective-planar graphs without a $K_{3,3}$-subdivision", submitted, arXiv:math.CO/\\0406140, (2004)].
\end{abstract}

\section{Introduction}
\noindent We use basic graph-theoretic terminology from Bondy and Murty \cite{Bondy} and Diestel \cite{Diestel}, and deal with undirected simple graphs.
Graph embeddings on a surface are important in VLSI design and in statistical mechanics. We are interested in non-planar graphs that can be embedded on the torus or on the projective plane.
By Kuratowski's theorem \cite{Kuratowski}, a graph $G$ is non-planar if and only if it contains a subdivision of $K_5$ or $K_{3,3}$ (see Figure~1). In this paper we characterize (and enumerate) the $2$-connected toroidal graphs with no $K_{3,3}$-subdivisions, following an analogous work for projective-planar graphs (\cite{GLL}). The next step in this research would be to characterize toroidal and projective-planar graphs containing a $K_{3,3}$-subdivision (with or without a $K_5$-subdivision).
\begin{figure}[h!] 
	\centerline {\includegraphics[width=2.8in]{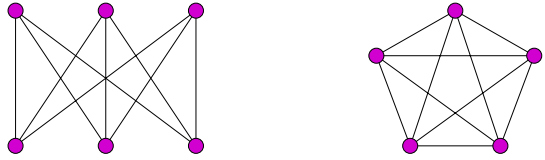}}
	\caption{Minimal non-planar graphs $K_{3,3}$ and $K_5$.}
\end{figure}

We assume that $G$ is a $2$-connected non-planar graph. A graph containing no $K_{3,3}$-subdivisions will be called {\it $K_{3,3}$-subdivision-free}. A general recursive decomposition of non-planar $K_{3,3}$-subdivision-free graphs is described in \cite{Wagner} and \cite{Kelmans}. A local decomposition of non-planar graphs containing a $K_5$-subdivision of a special type is described in \cite{Fellows} and \cite{GK} (some $K_{3,3}$-subdivisions are allowed), that is used later in \cite{GK} to detect a projective-planar or toroidal graph. The results of \cite{GK} provide a toroidality criterion for graphs containing a given $K_5$-subdivision and avoiding certain $K_{3,3}$-subdivisions by examining the embeddings of $K_5$ on the torus. The {\it torus} is an orientable surface of genus one which can be represented as a rectangle with two pairs of opposite sides identified. The graph $K_5$ has six different embeddings on the torus shown in Figure~2. Notice that the hatched region of each of the embeddings $E_1$ and $E_2$ forms a single face $F$.
\begin{figure}[h] 
	\centerline {\includegraphics[width=5.5in]{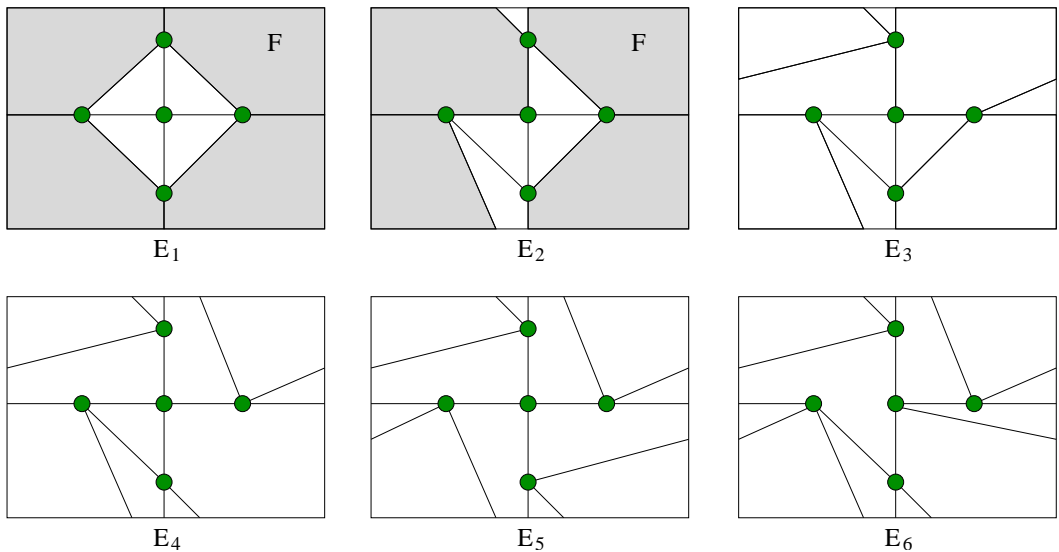}}
	\caption{Embeddings of $K_5$ on the torus.}
\end{figure}

In \cite{GLL} we prove the uniqueness of the decomposition of \cite{GK} for 2-connected non-planar {\em projective-planar} graphs with no $K_{3,3}$-subdivisions that gives a characterization of these graphs. In the present paper we state and prove an analogous structure theorem for the class $\T$ of 2-connected non-planar {\em toroidal} graphs with no $K_{3,3}$-subdivisions, involving certain ``circular crowns" of $K_5\backslash e$ networks and substitution of strongly planar networks for edges. The structure theorem provides a practical algorithm to recognize the toroidal graphs with no $K_{3,3}$-subdivisions in linear-time. Here we use the structure theorem to enumerate the labelled graphs in $\T$ by using the counting techniques of \cite{GLL} and \cite{Timothy} and improve known bounds for their number of edges. Finally, we enumerate the labelled graphs in $\T$ having no vertex of degree two. Tables can be found at the end of the paper.

\section{The structure theorem}
\noindent A {\it network} is a connected graph $N$ with two distinguished vertices $a$ and $b$, such that the graph $N\cup ab$ is $2$-connected. The vertices $a$ and $b$ are called the {\it poles} of $N$. The vertices of a network that are not poles are called {\it internal}. A network $N$ is {\it strongly planar} if the graph $N\cup ab$ is planar. We denote by $\N_P$ the class of strongly planar networks. 

The \emph{substitution} of a network $N$ for an edge $e=uv$ is done in the following way: 
choose an arbitrary orientation, say $\vec{e}=\vec{uv}$ of the edge, 
identify the pole $a$ of $N$ with the vertex $u$ and $b$ with $v$, 
and disregard the orientation of $e$ and the poles $a$ and $b$.
Note that both orientations of $e$ should be considered. 
It is assumed that the underlying set of $N$ is disjoint from $\{u,v\}$. 
The set of one or two resulting graphs is denoted by $e\uparrow N$.
More generally, given a graph $G_0$ with $k$ edges, $E=\{e_1, e_2, \ldots , e_k\}$, 
and a sequence $(N_1, N_2, \ldots , N_k)$ of disjoint networks, we define the {\it composition} $G_0\uparrow (N_1, N_2, \ldots , N_k)$ 
as the set of graphs that can be obtained by substituting the 
network $N_j$ for the edge $e_j$ of $G_0$, $j=1,2,\ldots ,k$. 
The graph $G_0$ is called the {\it core}, and the $N_i$'s are called the {\it components} 
of the resulting graphs. For a class of graphs $\G$ and a class of networks $\N$, we denote by $\G\uparrow \N$ the class of graphs obtained as compositions $G_0\uparrow (N_1, N_2, \ldots , N_k)$ with $G_0\in \G$ and $N_i\in \N$, $i=1,2,\ldots ,k$. 
We say that the composition $\G\uparrow \N$ is {\it canonical} if for any graph $G\in \G\uparrow \N$, 
there is a unique core $G_0\in \G$ and unique (up to orientation) components  
$N_1, N_2, \ldots , N_k\in \N$ that yield $G$.

In \cite{GLL} we prove the uniqueness of the representation $K_5\uparrow \N_P$ for $K_{3,3}$-subdivision-free projective-planar graphs. This gives an example of a canonical composition.
\begin{thm}[\cite{GK,GLL}]  A $2$-connected non-planar graph $G$ without a $K_{3,3}$-subdivision is projective-planar if and only if $G\in K_5\uparrow \N_P$. Moreover, the composition $K_5\uparrow \N_P$ is canonical.
\end{thm}
\begin{defn} {\em Given two $K_5$-graphs, the graph obtained by identifying an edge of one of the $K_5$'s with an edge of the other is called an {\it $M$-graph} (see Figure~3a)), and, when the edge of identification is deleted, an {\it $M^*$-graph} (see Figure~3b)).}
\end{defn}
\begin{figure}[h]
	\centerline {\includegraphics[width=3.0in]{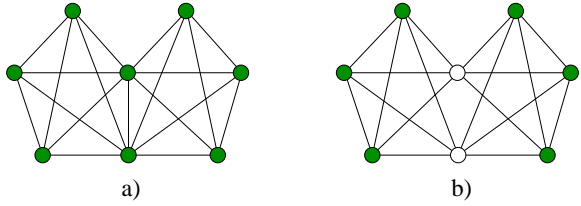}}
	\caption{a) $M$-graph, b) $M^*$-graph.}
\end{figure}
\begin{defn} {\em A network obtained from $K_5$ by removing the edge $ab$ between two poles is called a {\it $K_5\backslash e$-network}. A {\it circular crown} is a graph obtained from a cycle $C_i$, $i\ge 3$, by substituting $K_5\backslash e$-networks for some edges of $C_i$ in such  a way that no pair of  unsubstituted edges of $C_i$ are adjacent (see Figure~4).}
\end{defn}
\begin{figure}[h]
	\centerline {\includegraphics*[width=2.0in]{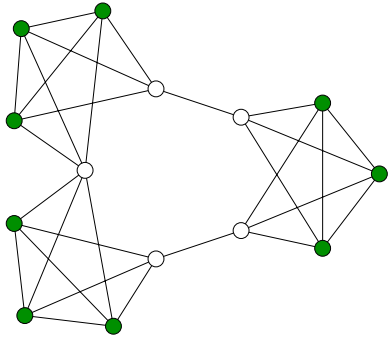}}
	\caption{A circular crown obtained from $C_5$.}
\end{figure}	
\begin{defn}
{\em A {\it toroidal core} is a graph $H$ which is isomorphic to either $K_5$, an $M$-graph, an $M^*$-graph, or a circular crown. We denote by $\T_C$ the class of toroidal cores.}
\end{defn}
The main result of this paper is the following structure theorem. The proof is given in Section 4.
\begin{thm}  A $2$-connected non-planar $K_{3,3}$-subdivision-free graph $G$ is toroidal if and only if $G\in \T_C\uparrow \N_P$. Moreover, the composition $\T=\T_C\uparrow \N_P$ is canonical.
\end{thm}
This theorem is used in Section 5 for the enumeration of labelled graphs in $\T$. In the future we hope to use Theorem 2 to enumerate unlabelled graphs in $\T$ as well.

\section{Related known results}
\noindent This section gives an overview of the structural results for toroidal graphs described in \cite{GK}. Following Diestel \cite{Diestel}, a $K_5$-subdivision is denoted by $TK_5$. The vertices of degree $4$ in $TK_5$ are the {\it corners\/} and the vertices of degree $2$ are the {\it inner vertices\/} of $TK_5$. For a pair of corners $a$ and $b$, the path $P_{ab}$ between $a$ and $b$  with all other vertices inner vertices is called a {\it side\/} of the $K_5$-subdivision. 

Let $G$ be a non-planar graph containing a fixed $K_5$-subdivision $TK_5$. A path $p$ in $G$ with one endpoint an inner vertex of $TK_5$, the other endpoint on a different side of
$TK_5$, and all other vertices and edges in $G\backslash TK_5$, is called a {\it short cut\/} of
the $K_5$-subdivision. A vertex $u\in G\backslash TK_5$ is called a
$3$-{\it corner vertex\/} with respect to $TK_5$ if $G\backslash TK_5$ contains internally disjoint paths connecting $u$ with at least three corners of the $K_5$-subdivision. 

\begin{propos} [\cite{Asano,Fellows,GK}] \label{propos:3-corner} 
Let $G$ be a non-planar graph  with a $K_5$-subdivision $TK_5$ for which there is either a short cut 
or a $3$-corner vertex. Then $G$ contains a $K_{3,3}$-subdivision.
\end{propos}
\begin{propos} [\cite{Fellows,GK}] \label{propos:sidecomponents}  
Let $G$ be a 2-connected graph with a $TK_5$ having no short cut or $3$-corner vertex.
Let $K$ denote the set of corners of $TK_5$.
Then any connected component $C$ of $G\backslash K$ contains inner vertices of at most one side
of $TK_5$ and $C$ is connected in $G$ to exactly two corners of $TK_5$. 
\end{propos}
Given a graph $G$ satisfying the hypothesis of Proposition \ref{propos:sidecomponents}, 
a \emph{side component} of $TK_5$ is defined as the subgraph of $G$ induced by a pair 
of corners $a$ and $b$ in $K$ and the connected components of $G\backslash K$ which are connected to both $a$ and $b$ in $G$. Notice that side components of $G$ can contain $K_{3,3}$-subdivisions. 
\begin{corol} [\cite{Fellows,GK}] \label{corol:unique corner}
For a 2-connected graph $G$ with a $TK_5$ having no short cut or $3$-corner vertex, 
two side components of $TK_5$ in $G$ have at most one vertex in common. 
The common vertex is the corner of intersection of two corresponding sides of $TK_5$.
\end{corol}
Thus we see that a graph $G$ satisfying the hypothesis of 
Proposition \ref{propos:sidecomponents}
can be decomposed into side components corresponding to the sides of $TK_5$.
Each side component $S$ contains exactly two corners $a$ and $b$ corresponding to a side of $TK_5$. 
If the edge $ab$ between the corners is not in $S$, we can add it to $S$ to obtain $S\cup ab$. 
Otherwise $S\cup ab=S$. We call $S\cup ab$ an \emph{augmented side component} of $TK_5$. Side components of a subdivision of an $M$-graph are defined by analogy with the side components of a $K_5$-subdivision by considering pairs of adjacent vertices of the $M$-graph. 

A planar side component $S$ of $TK_5$ in $G$ with two corners $a$ and $b$ is called {\it cylindrical} if the edge $ab\not\in S$ and the augmented side component $S\cup ab$ is non-planar. Notice that a planar side component $S=S\backslash ab$ is embeddable in a cylindrical section of the torus. A cylindrical section is provided by the face $F$ of the embeddings $E_1$ and $E_2$ of $K_5$ on the torus shown in  Figure~2. Toroidal graphs described in \cite{GK} can contain $K_{3,3}$-subdivisions because of a cylindrical side component $S$. An example of an embedding of the cylindrical side component $S=K_{3,3}\backslash e$ of a $TK_5$ on the torus is shown in Figure~6 where the graph $G$ of Figure~5 is embedded by completing the embedding $E_1$ of $K_5$ shown in Figure~2.
\begin{figure}[h!]
	\centerline {\includegraphics[width=1.9in]{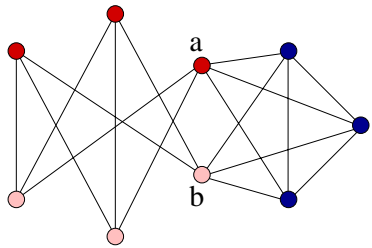}}
	\caption{A toroidal graph $G$ containing subdivisions of $K_{3,3}$ and of $K_5$.}
\end{figure}
\begin{figure}[h!]
	\centerline {\includegraphics{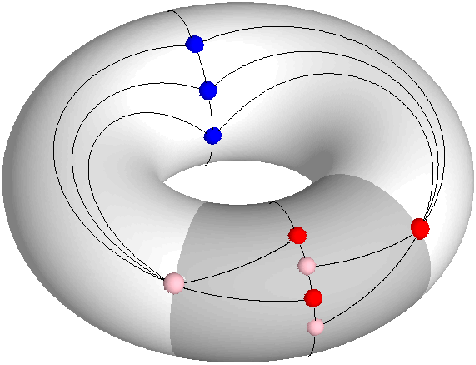}}
	\caption{Embedding of the cylindrical side component $K_{3,3}\backslash e$.}
\end{figure}

If a graph $G$ has no $K_{3,3}$-subdivisions, then Proposition \ref{propos:sidecomponents} 
can be applied, in virtue of Proposition \ref{propos:3-corner}. In this case, a result of \cite{GK} can be summarized as follows.
\begin{propos}[\cite{GK}] A 2-connected non-planar $K_{3,3}$-subdivision-free graph $G$ containing a $K_5$-sub\-division $TK_5$ is toroidal if and only if:

$(i)$ all the augmented side components of $TK_5$ in $G$ are planar graphs, or

$(ii)$ nine augmented side components of $TK_5$ in $G$ are planar, and  the remaining side component $S$ is cylindrical, or

$(iii)$ $G$ contains a subdivision $TM$ of an $M$-graph, and all the augmented side components of $TM$ in $G$ are planar.
\end{propos}
	
Further analysis of the cylindrical side component $S$ of Proposition 3(ii) will provide a proof of Theorem 2. Notice that graphs with $6$ or more vertices satisfying Propositon 3 are not $3$-connected. Therefore a $3$-connected non-planar graph different from $K_5$ must contain a $K_{3,3}$-subdivision (see also \cite{Asano}).

\section{Proof of the structure theorem}
\noindent 
A side component $S$ having two corners $a$ and $b$ can be considered as a network. We use the notation ${\rm Int}(S)$ to denote the {\em interior} of $S$, that is the subgraph ${\rm Int}(S)=S\backslash (\{a\}\cup \{b\})$ obtained by removing the two vertices $a$ and $b$. A network $S$ is called {\it cylindrical} if  $ab\not\in S$, $S$ is a planar graph, but $S\cup ab$ is non-planar. Recall that a network $S$ is called strongly planar if $S\cup ab$ is planar.

A {\em block} is a maximal 2-connected subgraph of a graph. A description of the {\em block-cutvertex tree decomposition} of a connected graph can be found in \cite{Diestel}. We consider blocks $G_i$ having two distinguished vertices $a_i$ and $b_i$. The distinguished vertices are called {\it poles} of the block. 
\begin{propos} Let $G$ be a 2-connected non-planar toroidal $K_{3,3}$-subdivision-free graph satisfying Proposition 3(ii) with the cylindrical side component $S$ having corners $a$ and $b$. Then the block-cutvertex decomposition of $S$ forms a path of blocks $S_1, S_2, \ldots ,S_k, k\ge 1$, as in Figure~7, and at least one of the blocks $S_1, S_2, \ldots ,S_k, k\ge 1$, is a cylindrical network. Moreover, every block $S_i$, $i=1,2,\ldots, k$, of $S$ is either a strongly planar network, or a cylindrical network of the form $K_5\backslash e\uparrow (N_1, N_2, \ldots, N_9)$, where $e=a_ib_i$ and the $N_j$'s are strongly planar networks.
\end{propos}
\proof Since $G$ is $2$-connected, each cut-vertex of $S$ belongs to exactly two blocks and lies on the corresponding side $P_{ab}$ of $TK_5$. Therefore the blocks of $S$ form a path as in Figure~7.
\begin{figure}[h!]
	\centerline {\includegraphics[width=4in]{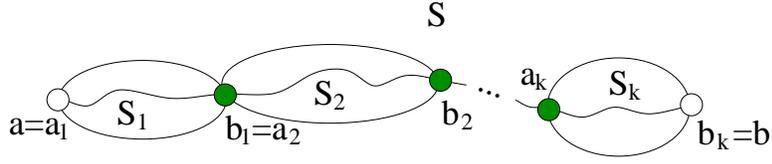}}
	\caption{\it Block-cutvertex decomposition for the cylindrical side component $S$.}
\end{figure}

Suppose each block $S_i$ of $S$, $i=1,2,\ldots ,k$, remains planar when the edge $a_ib_i$ is added to $S_i$. Then, clearly, $S\cup ab$ remains planar as well. Hence the fact that $S$ is cylindrical implies that at least one of the blocks $S_i$, $i=1,2,\ldots ,k$, is itself a cylindrical network.

Suppose a block $S_m$, $1\le m\le k$, of $S$ is cylindrical. Then, by Kuratowski's theorem, $S_m\cup a_mb_m$ contains a $K_5$-subdivision $TK'_5$. Clearly, $a_mb_m\in TK'_5$, $TK'_5$ has no short-cut or 3-corner vertex in $G$ and $a_m$ and $b_m$ are two corners of the $TK'_5$. The edge $a_mb_m$ of $TK'_5$ can be replaced by a path $P_{a_mb_m}$ in $G\backslash {\rm Int}(S_m)$ and we can decompose $G$ into the side components of $TK'_5$.

Since $G$ is toroidal and the side component $G\backslash {\rm Int}(S_m)$ of $TK'_5$ is cylindrical, all the other side components of $TK'_5$ in $G$ must be strongly planar networks by Proposition 3(ii). Therefore $S_m$ is a cylindrical network of the form $K_5\backslash e\uparrow (N_1, N_2, \ldots, N_9)$, with $e=a_mb_m$ and $N_j\in \N_P, j=1,2,\ldots,9$.
\hfill\rule{2mm}{2mm}

\bigbreak
Now we are ready to prove the structure Theorem 2 using Propositions 3 and 4.

\medbreak
\noindent {\bf Proof of Theorem 2.} (Sufficiency) Suppose $G$ is a graph in $\T_C\uparrow \N_P$, i.e. $G=H\uparrow (N_1, N_2, \ldots, N_k)$, where $H$ is a toroidal core having $k$ edges and $N_i$'s, $i=1,2,\ldots k$, are  strongly planar networks. If $H=K_5$ or $H=M$, then $G$ can be decomposed into the side components of $TK_5$ or $TM$ respectively and the augmented side components are planar graphs. Therefore, by Proposition 3(i) or 3(iii) respectively, $G$ is toroidal $K_{3,3}$-subdivision-free. 

If $H=M^*$ or $H$ is a circular crown, then we can choose a $K_5\backslash e$-network $N$ in $H$ and find a path $P_{ab}$ connecting $a$ and $b$ in the complementary part $H\backslash {\rm Int}(N)$. This determines a subdivision $TK_5$ in $G$ such that nine augmented side components of $TK_5$ in $G$ are planar, and the remaining side component $S$ defined by the corners $a$ and $b$ of $TK_5$ is cylindrical. Therefore, by Proposition 3(ii), $G$ is toroidal $K_{3,3}$-subdivision-free.
\bigbreak

(Necessity and Uniqueness) Let $G$ be a $2$-connected non-planar $K_{3,3}$-subdivision-free toroidal graph $G$. By Kuratowski's theorem, $G$ contains a $K_5$-subdivision $TK_5$. Let us prove that $G\in \T_C\uparrow \N_P$ by using Propositions 3 and 4. The fact that the composition $H\uparrow \N_P$, $H\in \T_C$, of $G$ is canonical will follow from the uniqueness of the sets of corner vertices in Proposition 3.

 Clearly, the sets of graphs corresponding to the cases (i), (ii) and (iii) of Proposition 3 are mutually disjoint. Suppose $G$ contains a subdivision $TK_5$ or $TM$ and all the augmented side components of $TK_5$ or $TM$, respectively, in $G$ are planar graphs as in Proposition 3(i, iii). Then $G=K_5\uparrow (N_1,N_2,\ldots ,N_{10})$ or $G=M\uparrow (N_1,N_2,\ldots ,N_{19})$, respectively, $K_5, M\in \T_C$ and all the $N_j$'s are in $\N_P$. The uniqueness of the decomposition in cases (i) and (iii) of Proposition 3 can be proved by analogy with Theorem 3 in \cite{GLL}: the set of corners of the $K_5$-subdivision in Proposition 3(i) and the set of corners of the $M$-graph subdivision in Proposition 3(iii) are uniquely defined. This covers toroidal cores $K_5$ and the $M$-graph.

Suppose $S$ is the unique cylindrical side component of $TK_5$ in $G$ as in Proposition 3(ii). Notice that $G\backslash {\rm Int}(S)$ itself is a cylindrical network of the form $K_5\backslash e\uparrow (N_1, N_2, \ldots, N_9)$, where $e=ab$ and $N_j\in \N_P, j=1,2,\ldots,9$. By Proposition 4, the block-cutvertex decomposition of $S$ forms a path of blocks $S_1, S_2, \ldots ,S_k, k\ge 1$, as in Figure~7, and at least one of the blocks $S_1, S_2, \ldots ,S_k, k\ge 1$, is a cylindrical network. In this path we can regroup maximal series of consecutive strongly planar networks into single strongly planar networks so that at most one strongly planar network $N^\prime$ is separating two cylindrical networks in the resulting path, and the poles of the strongly planar network $N^\prime$ are uniquely defined by maximality. By Proposition 4, the cylindrical networks in the path are of the form $K_5\backslash e\uparrow (N_1, N_2, \ldots, N_9)$, where $N_j\in \N_P, j=1,2,\ldots,9$, and the corners $a'$ and $b'$, $e=a'b'$, are uniquely defined with respect to the corresponding $K_5$-subdivision $TK'_5$ in $G$. Therefore the unique set of corners completely defines a toroidal core $M^*$ or a circular crown $H$ having $k$ edges and the set of corresponding strongly planar networks $N_1,N_2,\ldots ,N_k$, such that $G=M^*\uparrow (N_1,N_2,\ldots ,N_{18})$ or $G=H\uparrow (N_1,N_2,\ldots ,N_k)$, respectively.
\hfill\rule{2mm}{2mm}\\

Theorems 1 and 2 imply that a projective-planar graph with no $K_{3,3}$-subdivisions is toroidal. However an arbitrary projective-planar graph can be non-toroidal. The characterizations of Theorems 1 and 2 can be used to detect projective-planar or toroidal graphs with no $K_{3,3}$-subdivisions in linear time. The implementation of this algorithm can be derived from \cite{GK} by using a breadth-first or depth-first search technique for the decomposition and by doing a linear-time planarity testing. The linear-time complexity follows from the linear-time complexity of the decomposition and from the fact that each vertex of the initial graph can appear in at most 7 different components.

A corollary to Euler's formula for the plane says that a planar graph with $n\ge 3$ vertices 
can have at most $3n-6$ edges (see, for example, \cite{Bondy} and \cite{Diestel}). Let us state this for 2-connected planar graphs with $n$ vertices and $m$ edges as follows:
\begin{equation} \label{eq:mplanar}
m\le  \left\{ \begin{array}{ll} 
 3n-5 & \mbox{if \,$n=2$} \\ 
 3n-6 & \mbox{if \,$n\ge 3$} 
\end{array}.
\right. 
\end{equation} 
In fact, $m=3n-5=1$ if $n=2$. The generalized Euler formula (see, for example, \cite{Carsten}) implies 
that a toroidal graph $G$ with $n$ vertices can have up to $3n$ edges. 
An arbitrary graph $G$ without a $K_{3,3}$-subdivision is known to have at most $3n-5$ edges 
(see \cite{Asano}). The following proposition shows that toroidal graphs with no $K_{3,3}$-subdivisions satisfy a stronger relation, which is analogous to planar graphs. 
\begin{propos}  \label{propos:mtor}
The number $m$ of edges of a non-planar $K_{3,3}$-subdivision-free toroidal $n$-vertex graph $G$ 
satisfies $m\le 3n-5$ if $n=5 \ \mbox{or}\ 8$, and 
\begin{equation} \label{eq:mtor}
m\le 3n-6, \ \mbox{if}\  n\ge 6 \ \mbox{and}\ n\not=8.
\end{equation}
\end{propos}
\proof Clearly, toroidal graphs satisfying Theorem 2 also satisfy Proposition 3. By Proposition $3(i,ii)$, each side component $S_i$ of $TK_5$ in $G$, $i=1,2,\ldots,10$, 
satisfies the condition $(\ref{eq:mplanar})$ with $n=n_i$, the number of vertices, and $m=m_i$, the number of edges of $S_i$, $i=1,2,\ldots,10$. Since each corner of $TK_5$ is in precisely $4$ side components, we have $\sum_{i=1}^{10}n_i=n+ 15$ and we obtain, by summing these 10 inequalities, 
$$
m=\sum_{i=1}^{10}m_i\le \left\{ \begin{array}{lll} 
3\sum_{i=1}^{10}n_i-50 = 3(n+15)-50 = 3n-5 & \mbox{if $n=5$} \\ \\
3\sum_{i=1}^{10}n_i-51 = 3(n+15)-51 = 3n-6 & \mbox{if $n\ge 6$} 
\end{array},
\right.
$$
since 
$n=5$ iff $n_i=2$, $i=1,2,\ldots,10$, and $n\ge 6$ if and only if at least one $n_j\ge 3$, $j=1,2,\ldots,10$.

Similarly, by Proposition $3(iii)$, each side component $S_i$ of $TM$ in $G$, $i=1,2,\ldots,19$, 
satisfies the condition $(\ref{eq:mplanar})$ with $n=n_i$, the number of vertices, and $m=m_i$, the number of edges of $S_i$, $i=1,2,\ldots,19$. Since 2 vertices of $TM$ are in precisely $7$ side components, 6 vertices of $TM$ are in precisely $4$ side components, and all the other vertices of $G$ are in a unique side component, we have $\sum_{i=1}^{19}n_i=n+ 30$ and we obtain, by summing these 19 inequalities, 
$$
m=\sum_{i=1}^{19}m_i\le \left\{ \begin{array}{lll} 
3\sum_{i=1}^{19}n_i-95 = 3(n+30)-95 = 3n-5 & \mbox{if $n=8$} \\ \\
3\sum_{i=1}^{19}n_i-96 = 3(n+30)-96 = 3n-6 & \mbox{if $n\ge 9$} 
\end{array},
\right.
$$
since 
$n=8$ iff $n_i=2$, $i=1,2,\ldots,19$, and $n\ge 9$ if and only if at least one $n_j\ge 3$, $j=1,2,\ldots,19$.
\hfill\rule{2mm}{2mm}
\bigbreak
An analogous result for the projective-planar graphs can be found in \cite{GLL}. Also note that Corollary 8.3.5 of \cite{Diestel} implies that graphs 
with no $K_5$-minors can have at most $3n-6$ edges.

\section{Counting labelled $K_{3,3}$-sub\-division-free toroidal\\ graphs}
\noindent Now let us consider the question of the labelled enumeration of toroidal 
graphs with no $K_{3,3}$-subdivisions according to the numbers of vertices and edges.
First, we review some basic notions and terminology of labelled enumeration together with the counting methods and technique used in \cite{Timothy, GLL}. The reader should have some familiarity with exponential generating functions and their operations (addition, multiplication and composition).
For example, see \cite{BLL}, \cite{Goulden}, \cite{Stanley}, or \cite{Wilf}.

By a \emph{labelled} graph, we mean a simple graph $G=(V,E)$ where the set of vertices $V=V(G)$ 
is itself the set of labels and the labelling function is the identity function.
$V$ is called the \emph{underlying set} of $G$.
An edge $e$ of $G$ then consists of an unordered pair $e=uv$ of elements of $V$ 
and $E=E(G)$ denotes the set of edges of $G$.
If $W$ is another set and $\sigma:V\tilde{\rightarrow}W$ is a bijection, 
then any graph $G=(V,E)$ on $V$,  
can be transformed into a graph $G^\prime=\sigma(G)=(W,\sigma(E))$, where
$\sigma(E)=\{\sigma(e) =\sigma(u)\sigma(v)\,|\, e\in E\}$. 
We say that $G^\prime$ is obtained from $G$ by \emph{vertex relabelling}
and that $\sigma$ is a graph \emph{isomorphism} $G\tilde{\rightarrow}G^\prime$.
An \emph{unlabelled graph} is then seen as an isomorphism class $\gamma$ of labelled graphs. 
We write $\gamma=\gamma(G)$ if $\gamma$ is the isomorphism class of $G$. 
By the \emph{number of ways to label} an unlabelled graph $\gamma(G)$, 
where $G=(V,E)$, we mean the number of distinct graphs $G^\prime$ on 
the underlying set $V$ which are isomorphic to $G$. 
Recall that this number is given by $n!/|\Aut(G)|$, 
where $n=|V|$ and $\Aut(G)$ denotes the automorphism
group of $G$.

A \emph{species} of graphs is a class of labelled graphs which is closed under vertex relabellings. 
Thus any class $\G$ of unlabelled graphs gives rise to a species, also denoted by $\G$,
by taking the set union of the isomorphism classes in $\G$.
For any species $\G$ of graphs, we introduce its (exponential) \emph{generating function}
$\G(x,y)$ as the formal power series
\begin{equation} \label{eq:Gdexy}
\G(x,y)=\sum_{n\ge 0}g_n(y)\frac{x^n}{n!}, \ \ 
\mbox {with}\ \ g_n(y)=\sum_{m\geq0} g_{n,m}y^m,
\end{equation}
where $g_{n,m}$ is the number of graphs in $\G$ with $m$ edges over a given set
of vertices $V_n$ of size $n$. Here $y$ is a formal variable which acts as an edge counter.
For example, for the species $\G=K=\{K_n\}_{n\geq0}$ of complete graphs, we have
\begin{equation} \label{eq:Kxy}
K(x,y)=\sum_{n\geq0}y^{n\choose2}x^n/n!,
\end{equation}
while for the species $\G=\G_a$ of all simple graphs,
we have $\G_a(x,y)=K(x,1+y)$. 

A species of graphs is \emph{molecular} if it contains only one isomorphism class.
For a molecular species $\gamma=\gamma(G)$, where $G$ has $n$ vertices and $m$ edges, 
we have $\gamma(x,y)= \frac{y^mn!}{|\Aut(G)|}x^n/n!=y^mx^n/|\Aut(G)|$. For example,
\begin{eqnarray} \label{eq:K5dexy}
K_5(x,y)=\frac{x^5y^{10}}{5!}.
\end{eqnarray} 
Also, for the graphs $M$ and $M^*$ described in Section 2, we have
\begin{eqnarray} \label{eq:Mdexy}
M(x,y)=280\frac{x^8y^{19}}{8!}, \ \ \ M^*(x,y)=280\frac{x^8y^{18}}{8!},
\end{eqnarray} 
since $|\Aut(M)|=|\Aut(M^*)|=144$.

For the enumeration of networks, we consider that the poles $a$ and $b$ are not labelled, or, in other words, that only the internal vertices 
form the underlying set. Hence the generating function of a class (or species) $\N$ of networks is defined by 
\begin{equation} \label{eq:Ndexy}
\N(x,y)=\sum_{n\ge 0}\nu_n(y)\frac{x^n}{n!}, \ \ 
\mbox {with}\ \ \nu_n(y)=\sum_{m\geq0} \nu_{n,m}y^m,
\end{equation}
where $\nu_{n,m}$ is the number of networks in $\N$ with $m$ edges and a given set
of internal vertices $V_n$ of size $n$. For example, we have 
\begin{equation} \label{eq:K5net}
(K_5\backslash e)(x,y)=\frac{x^3y^{9}}{3!},
\end{equation}

A species $\cal{N}$ of networks is called {\it symmetric} if for any $\cal{N}$-network 
$N$ (i.e. $N$ in $\cal{N}$), the \emph{opposite network} $\tau\cdot N$, obtained by interchanging 
the poles $a$ and $b$, is also in $\cal{N}$. Examples of symmetric species of networks are the classes $\N_P$, of strongly planar networks, and $\R$, of series-parallel networks (see \cite{Timothy, GLL}).

\begin{lemma} \emph{(T. Walsh \cite{Timothy, GLL})} \label{lemme:GflecheN}
Let $\G$ be a species of graphs and $\N$ be a symmetric species of networks 
such that the composition $\G\uparrow \N$ is canonical. 
Then the following generating function identity holds:
\begin{equation} \label{eq:GflecheNxy}
(\G\uparrow \N)(x,y)=\G(x,\N(x,y)).
\end{equation}
\end{lemma}
By Theorem 2 and Lemma 1, we have the following proposition. 
\begin{propos} 
The generating function $\T(x,y)$  of labelled non-planar $K_{3,3}$-subdivision-free toroidal graphs is given by
\begin{equation} \label{eq:Txy}
\T(x,y)=(\T_C\uparrow \N_P)(x,y)=\T_C(x,\N_P(x,y)),
\end{equation}
where $\T_C$ denotes the class of toroidal cores (see Definition 3).
\end{propos}

Let $P$ denote the species of $2$-connected planar graphs. Then the generating function of $\N_P$, the associated class of strongly planar networks, is given by
\begin{equation} \label{eq:NPxy}
\N_P(x,y) = (1+y)\frac{2}{x^2}\frac{\partial}{\partial y}P(x,y) - 1
\end{equation}
(see \cite{Timothy, GLL}). Methods for computing the generating function $P(x,y)$ of labelled $2$-connected planar graphs
are described in \cite{Bender} and \cite{Bodirsky}. Formula (\ref{eq:NPxy}) can then be used to compute $\N_P(x,y)$. Therefore there remains only to compute the generating function $\T_C(x,y)$ for toroidal cores. Recall that $\T_C=K_5+M+M^*+CC$, where $CC$ denotes the class of circular crowns. Circular crowns can be enumerated as follows using matching polynomials.
\begin{propos}
The mixed generating series $CC(x,y)$ of circular crowns is given by 
\begin{eqnarray} \label{eq:ccrowns}
CC(x,y)=-\frac{12x^4y^9+12x^5y^{10}+x^8y^{18}+72\ln(1-\frac{x^4y^9}{6}-\frac{x^5y^{10}}{6})}{144}.
\end{eqnarray}
\end{propos}
\proof Recall that a {\em matching} $\mu$ of a finite graph $G$ is a set of disjoint edges of $G$. We define the {\em matching polynomial} of $G$ as
\begin{eqnarray}
M_G(y)=\sum_{\mu \in \M(G)}y^{|\mu|},
\end{eqnarray}
where $\M(G)$ denotes the set of matchings of $G$. In particular, the matching polynomials $U_n(y)$ and $T_n(y)$ for paths and cycles of size $n$ are well known (see \cite{Godsil}). They are closely related to the Chebyshev polynomials. To be precise, let $P_n$ denote the path graph $(V,E)$ with $V=[n]=\{1,2,\ldots,n\}$ and $E=\{\{i,i+1\} |\ i=1,2,\ldots,n-1\}$ and $C_n$ denote the cycle graph with $V=[n]$ and $E=\{\{i,i+1(mod\ n)\} |\ i=1,2,\ldots,n\}$. Then we have 
\begin{equation} \label{eq:Undey}
U_n(y)=\sum_{\mu \in \M(P_n)}y^{|\mu|},\ \ \ T_n(y)=\sum_{\mu \in \M(C_n)}y^{|\mu|}.
\end{equation}
The dichotomy caused by the membership of the edge $\{n-1,n\}$ in the matchings of the path $P_n$ leads to the recurrence relation
\begin{equation} \label{eq:Undeyrec}
U_n(y)=y U_{n-2}(y)+U_{n-1}(y),
\end{equation}
for $n\ge 2$, with $U_0(y)=U_1(y)=1$. It follows that the ordinary generating function of the matching polynomials $U_n(y)$ is rational. In fact, it is easily seen that  
\begin{equation} \label{eq:sumUndey}
\sum_{n\ge 0}U_n(y)x^n=\frac{1}{1-x-yx^2}.
\end{equation}

Now, the dichotomy caused by the membership of the edge $\{1,n\}$ in the matchings of the cycle $C_n$ leads to the relation
\begin{equation} \label{eq:Tndeyrec}
T_n(y)=y U_{n-2}(y)+U_{n}(y),
\end{equation}
for $n\ge 3$. It is then a simple matter, using (\ref{eq:sumUndey}) and (\ref{eq:Tndeyrec}) to compute their ordinary generating function, denoted by $G(x,y)$. We find 
\begin{equation} \label{eq:sumTndey}
G(x,y)=\sum_{n\ge 3}T_n(y)x^n=\frac{x^3(1+3y+yx+2y^2x)}{1-x-yx^2}.
\end{equation}
In fact, we also need to consider the {\em homogeneous matchings polynomials}
\begin{equation} \label{eq:TwoVariables}
T_n(y,z)=z^nT_n(\frac{y}{z})=\sum_{\mu \in \M(C_n)}y^{|\mu|}z^{n-|\mu|},
\end{equation}
where the variable $z$ marks the edges which are not selected in the matchings, whose generating function $G(x,y,z)=\sum_{n\ge 3}T_n(y,z)x^n$ is given by
\begin{equation} \label{eq:Gxyz}
G(x,y,z)=G(xz,\frac{y}{z})=\frac{x^3z^2(z+3y+xyz+2xy^2)}{1-xz-x^2yz}.
\end{equation}

We now introduce the species $BC$ of pairs $(c,\mu)$, where $c$ is a cycle of length $n\ge 3$ and $\mu$ is a matching of $c$, with weight $y^{|\mu|}z^{n-|\mu|}$. Since there are $\frac{(n-1)!}{2}$ non-oriented cycles on a set of size $n\ge 3$, and all these cycles admit the same homogeneous matching polynomial $T_n(y,z)$, the mixed generating function of labelled $BC$-structures is 
\begin{eqnarray} \label{eq:bicycles}
BC(x,y,z)=\sum_{n\ge 3}\frac{(n-1)!}{2}T_n(y,z)\frac{x^n}{n!} \nonumber \\
=\frac{1}{2}\sum_{n\ge 3}T_n(y,z)\frac{x^n}{n} \nonumber \\
=\frac{1}{2}\int_{0}^{x}\frac{1}{t}G(t,y,z)\,dt
 \nonumber \\
=-\frac{2xz+2x^2zy+x^2z^2+2\ln(1-xz-x^2yz)}{4}.
\end{eqnarray}

Notice that in a circular crown, the unsubstituted edges are not adjacent, by definition, and hence form a matching of the underlying cycle, while the substituted edges are replaced by $K_5\backslash e$-networks. We can thus write 
\begin{eqnarray} \label{eq:DecCircCr}
CC=BC\uparrow_z(K_5\backslash e),
\end{eqnarray}
where the notation $\uparrow_z$ means that only the edges marked by $z$ are replaced by $K_5\backslash e$-networks. Moreover the decomposition (\ref{eq:DecCircCr}) is canonical and we have
\begin{eqnarray} \label{eq:CircCr}
CC(x,y)=BC(x,y,(K_5\backslash e)(x,y)),
\end{eqnarray}
which implies (\ref{eq:ccrowns}) using (\ref{eq:K5net}).
\hfill\rule{2mm}{2mm}\\

A substitution of the generating function $\N_P(x,y)$ of (\ref{eq:NPxy}) counting the strongly planar networks for the variable $y$ in (\ref{eq:Mdexy}), (\ref{eq:K5dexy}), and (\ref{eq:ccrowns}) gives the generating function for labelled $2$-connected non-planar toroidal graphs with no $K_{3,3}$-subdivision, i.e.
\begin{equation} \label{eq:TorGr}
\T(x,y)=K_5(x,\N_P(x,y))+M(x,\N_P(x,y))+M^*(x,\N_P(x,y))+CC(x,\N_P(x,y)).
\end{equation}

Notice that the term $K_5(x,\N_P(x,y))$ in (\ref{eq:TorGr}) also enumerates non-planar 2-connected $K_{3,3}$-subdivision-free {\em projective-planar graphs} and that corresponding tables are given in \cite{GLL}. Here we present the computational results just for labelled graphs in $\T$ that are not projective-planar.  Numerical results are presented in Tables~1 and~2, where $\T(x,y)-K_5(x,\N_P(x,y))=\sum_{n\geq8}\sum_{m} t_{n,m}x^ny^m/n!$ and $t_n=\sum_{m}t_{n,m}$ count labelled non-projective-planar graphs in $\T$.

The homeomorphically irreducible non-projective-planar graphs in $\T$, i.e. the graphs having no vertex of degree two, can be counted by using several methods described in detail in Section 4 of \cite{GLL}. We used the approach of Proposition 8 of \cite{GLL} to obtain the numerical data presented in Tables~3 and~4 for labelled homeomorphically irreducible graphs in $\T$ that are not projective-planar.

%
\begin{table}[!h] \label{table:fnm}
\centerline {\footnotesize
	\begin{tabular}[t]{|| r | r | r || r | r | r || r | r | r ||}
	\hline
	$n$ & $m$ & $t_{n,m}$ & $n$ & $m$ & $t_{n,m}$ & $n$ & $m$ & $t_{n,m}$\\
	\hline \hline
	8 & 18 & 280 &  13 & 23 & 1838008972800 & 15 & 25 & 5973529161600000\\
	\hline
	8 & 19 & 280 & 13 & 24 & 12383684913600 & 15 & 26 & 60679359861120000\\
	\hline
	9 & 19 & 50400 & 13 & 25 & 36576568828800 & 15 & 27 & 280619124786000000\\
	\hline
	9 & 20 & 93240 & 13 & 26 & 61986597472800 & 15 & 28 & 785755439324856000\\
	\hline
	9 & 21 & 47880 & 13 & 27 & 66199273620480 & 15 & 29 & 1496142328612932000\\
	\hline
	10 & 20 & 5292000 & 13 & 28 & 46419992138520 & 15 & 30 & 2068477720590481200\\
	\hline
	10 & 21 & 15044400 & 13 & 29 & 22180672954440 & 15 & 31 & 2175937397296462800\\
	\hline
	10 & 22 & 15510600 & 13 & 30 & 7737403073400 & 15 & 32 & 1810128996903427200\\
	\hline
	10 & 23 & 5972400 & 13 & 31 & 2053743892200 & 15 & 33 & 1223242124356652400\\
	\hline
	10 & 24 & 239400 & 13 & 32 & 348540192000 & 15 & 34 & 673154380612513800\\
	\hline
	11 & 21 & 426888000 & 13 & 33 & 27935107200 & 15 & 35 & 293316332440131000\\
	\hline
	11 & 22 & 1700899200 & 14 & 24 & 107217190080000 & 15 & 36 & 96295664217753000\\
	\hline
	11 & 23 & 2724044400 & 14 & 25 & 896474952172800 & 15 & 37 & 22260497063805000\\
	\hline
	11 & 24 & 2136842400 & 14 & 26 & 3359265613704000 & 15 & 38 & 3218036781960000\\
	\hline
	11 & 25 & 773295600 & 14 & 27 & 7460402644094400 & 15 & 39 & 218263565520000\\
	\hline
	11 & 26 & 94386600 & 14 & 28 & 10948159170748800 & 16 & 26 & 322570574726400000\\
	\hline
	11 & 27 & 7900200 & 14 & 29 & 11253868616390400 & 16 & 27 & 3914073525922560000\\
	\hline
	12 & 22 & 29455272000 & 14 & 30 & 8467602606022560 & 16 & 28 & 21877169871997440000\\
	\hline
	12 & 23 & 155542464000 & 14 & 31 & 4876995169606560 & 16 & 29 & 75157668529175232000\\
	\hline
	12 & 24 & 348414066000 & 14 & 32 & 2222245323698400 & 16 & 30 & 178928606393593056000\\
	\hline
	12 & 25 & 424294516800 & 14 & 33 & 785187373370400 & 16 & 31 & 316283670286218835200\\
	\hline
	12 & 26 & 297599563800 & 14 & 34 & 197208318106800 & 16 & 32 & 435483254883942064800\\
	\hline
	12 & 27 & 118905448200 & 14 & 35 & 31064455422000 & 16 & 33 & 484253520685973438400\\
	\hline
	12 & 28 & 27683548200 & 14 & 36 & 2294786894400 & 16 & 34 & 445576710488584474800\\
	\hline
	12 & 29 & 4821201000 & \multicolumn{3}{c ||}{ } & 16 & 35 & 341998556200139638800\\
	\cline {1-3}\cline {7-9}
	12 & 30 & 410810400 & \multicolumn{3}{c ||}{ } & 16 & 36 & 216864722075241240000 \\
	\cline {1-3}\cline {7-9}
	\multicolumn{6}{c ||}{ } & 16 & 37 & 111029372376938215200\\
	\cline {7-9}
	\multicolumn{6}{c ||}{ } & 16 & 38 & 44479356838490574000\\
	\cline {7-9}
	\multicolumn{6}{c ||}{ } & 16 & 39 & 13374653821603074000\\
	\cline {7-9}
	\multicolumn{6}{c ||}{ } & 16 & 40 & 2831094029443680000\\
	\cline {7-9}
	\multicolumn{6}{c ||}{ } & 16 & 41 & 375386906774880000\\
	\cline {7-9}
	\multicolumn{6}{c ||}{ } & 16 & 42 & 23417178744960000\\
	\cline {7-9}
	\end{tabular} }
	\caption{The number of labelled non-planar non-projective-planar toroidal 2-connected graphs without a $K_{3,3}$-sub\-di\-vi\-sion (having $n$ vertices and $m$ edges).}
\end{table}
%
\begin{table}[!h] \label{table:fn}
\centerline {\footnotesize
	\begin{tabular}[t]{|| r | r ||}
	\hline
	$n$ & $t_n$ \\
	\hline \hline
	8 & 560 \\
	\hline
	9 & 191520\\
	\hline
	10 & 42058800\\
	\hline
	11 & 7864256400\\
	\hline
	12 & 1407126890400\\
	\hline
	13 & 257752421166240\\
	\hline
	14 & 50607986220311520\\
	\hline
	15 & 10995419195575214400\\
	\hline
	16 & 2692773804667509763200\\
	\hline
	17 & 747221542837742897724800\\
	\hline
	18 & 233698171655650029030743040\\
	\hline
	19 & 81472765051132560093387934080\\
	\hline
	20 & 31268587126068905034073041062400\\
	\hline
	\end{tabular} }
	\caption{The number of labelled non-planar non-projective-planar toroidal 2-connected $K_{3,3}$-subdivision-free graphs (having $n$ vertices).}
\end{table}

%
\begin{table}[!h] \label{table:fnm}
\centerline {\footnotesize
	\begin{tabular}[t]{|| r | r | r || r | r | r || r | r | r ||}
	\hline
	$n$ & $m$ & $t_{n,m}$ & $n$ & $m$ & $t_{n,m}$ & $n$ & $m$ & $t_{n,m}$\\
	\hline \hline
	8 & 18 & 280 &  14 & 26 & 6054048000 & 16 & 29 & 5811886080000\\
	\hline
	8 & 19 & 280 & 14 & 27 & 285751065600 & 16 & 30 & 621544891968000\\
	\hline
	9 & 19 & 5040 & 14 & 28 & 3361812854400 & 16 & 31 & 11935943091072000\\
	\hline
	10 & 20 & 25200 & 14 & 29 & 17840270448000 & 16 & 32 & 101350194001056000\\
	\hline
	10 & 22 & 226800 & 14 & 30 & 55133382704400 & 16 & 33 & 499371733276416000\\
	\hline
	10 & 23 & 466200 & 14 & 31 & 108994658572800 & 16 & 34 & 1611221546830896000\\
	\hline
	10 & 24 & 239400 & 14 & 32 & 141179453415000 & 16 & 35 & 3605404135132800000\\
	\hline
	11 & 23 & 10256400 & 14 & 33 & 118498240060200 & 16 & 36 & 5738963267481444000\\
	\hline
	11 & 24 & 30492000 & 14 & 34 & 61801664324400 & 16 & 37 & 6540526990277280000\\
	\hline
	11 & 25 & 43520400 & 14 & 35 & 18158435895600 & 16 & 38 & 5293490794557966000\\
	\hline
	11 & 26 & 31185000 & 14 & 36 & 2294786894400 & 16 & 39 & 2967845927880834000\\
	\hline
	11 & 27 & 7900200 & 15 & 28 & 1961511552000 & 16 & 40 & 1095216458944608000\\
	\hline
	12 & 24 & 189604800 & 15 & 29 & 57537672192000 & 16 & 41 & 239190441890400000\\
	\hline
	12 & 25 & 1079416800 & 15 & 30 & 557188343712000 & 16 & 42 & 23417178744960000\\
	\hline
	12 & 26 & 3044487600 & 15 & 31 & 2827950253128000 & 17 & 31 & 3903916528512000\\
	\hline
	12 & 27 & 5080614000 & 15 & 32 & 8936155496268000 & 17 & 32 & 174648084811200000\\
	\hline
	12 & 28 & 4776294600 & 15 & 33 & 18886100303070000 & 17 & 33 & 2606052624215040000\\
	\hline
	12 & 29 & 2261536200 & 15 & 34 & 27395286118200000 & 17 & 34 & 20178959825344320000\\
	\hline
	12 & 30 & 410810400 & 15 & 35 & 27296971027326000 & 17 & 35 & 97287841256493888000\\
	\hline
	13 & 25 & 1686484800 & 15 & 36 & 18324093378591000 & 17 & 36 & 319780940570307216000\\
	\hline
	13 & 26 & 22875652800 & 15 & 37 & 7906712877063000 & 17 & 37 & 751384930811218704000\\
	\hline
	13 & 27 & 126680954400 & 15 & 38 & 1978851858984000 & 17 & 38 & 1292496613555066920000\\
	\hline
	13 & 28 & 382608626400 & 15 & 39 & 218263565520000 & 17 & 39 & 1642597679422623924000\\
	\hline
	13 & 29 & 700723623600 & \multicolumn{3}{c ||}{ } & 17 & 40 & 1539140405659676820000\\
	\cline {1-3}\cline {7-9}
	13 & 30 & 788388400800 & \multicolumn{3}{c ||}{ } & 17 & 41 & 1049167407329489448000\\
	\cline {1-3}\cline {7-9}
	13 & 31 & 525156231600 & \multicolumn{3}{c ||}{ } & 17 & 42 & 505608857591934096000 \\
	\cline {1-3}\cline {7-9}
	13 & 32 & 188324136000 & \multicolumn{3}{c ||}{ } & 17 & 43 & 163183484418946992000\\
	\cline {1-3}\cline {7-9}
	13 & 33 & 27935107200 & \multicolumn{3}{c ||}{ } & 17 & 44 & 31635477128166912000\\
	\cline {1-3}\cline {7-9}
	\multicolumn{6}{c ||}{ } & 17 & 45 & 2784602773016064000\\
	\cline {7-9}
	\end{tabular} }
	\caption{The number of labelled non-planar non-projective-planar toroidal 2-connected $K_{3,3}$-subdivision-free graphs with no vertex of degree 2 (having $n$ vertices and $m$ edges).}
\end{table}
%
\begin{table}[!h] \label{table:fn}
\centerline {\footnotesize
	\begin{tabular}[t]{|| r | r ||}
	\hline
	$n$ & $t_n$ \\
	\hline \hline
	8 & 560 \\
	\hline
	9 & 5040\\
	\hline
	10 & 957600\\
	\hline
	11 & 123354000\\
	\hline
	12 & 16842764400\\
	\hline
	13 & 2764379217600\\
	\hline
	14 & 527554510282800\\
	\hline
	15 & 114387072405606000\\
	\hline
	16 & 27728561968887780000\\
	\hline
	17 & 7418031804967840056000\\
	\hline
	18 & 2167306256125914230527200\\
	\hline
	19 & 685709965521372865035362400\\
	\hline
	20 & 233306923207078035272369412000\\
	\hline
	\end{tabular} }
	\caption{The number of labelled non-planar non-projective-planar toroidal 2-connected $K_{3,3}$-subdivision-free graphs with no vertex of degree 2 (having $n$ vertices).}
\end{table}

\end{document}